\documentclass[11pt,a4paper]{article}
\usepackage{mathrsfs}
\usepackage{amsfonts}
\usepackage{}

\usepackage{graphicx}
\usepackage{amsfonts,amsmath, amssymb}

\newtheorem{theo}{\textbf{\ \ \quad Theorem}}[section]

\newtheorem{lem}{\textbf{\ \ \quad Lemma}}[section]
\newtheorem{remark}{\textbf{\ \ \quad Remark}}[section]

\newtheorem{prop}{\textbf{\ \ \quad Proposition}}[section]
\newtheorem{defi}{\textbf{\ \ \quad Definition}}[section]

\newcommand{\lbl}[1]{\label{#1}}

\newcommand{\be}{\begin{equation}}
\newcommand{\ee}{\end{equation}}
\newcommand\bes{\begin{eqnarray}}
\newcommand\ees{\end{eqnarray}}
\newcommand{\bess}{\begin{eqnarray*}}
\newcommand{\eess}{\end{eqnarray*}}

\newcommand{\nm}{\nonumber}

{\setlength\arraycolsep{2pt}}

\title{Regularity of stochastic nonlocal diffusion equations }

\author{Guangying Lv$^a$, Hongjun Gao$^b$ Jinlong Wei$^c$, Jiang-Lun Wu$^d$\\
\\
\ \\
   {\small \it $^a$Institute of Contemporary Mathematics, Henan University}\\
  {\small \it Kaifeng, Henan 475001, China}\\
  {\small \tt gylvmaths@henu.edu.cn}\\
    {\small \it $^b$Institute of Mathematics, School of Mathematical Science}\\
  {\small \it Nanjing Normal University, Nanjing 210023, China}\\
  {\small \tt gaohj@njnu.edu.cn}\\
  {\small \it $^c$ School of Statistics and Mathematics, Zhongnan University of}\\
  {\small \it
 Economics and Law, Wuhan, Hubei 430073, China}\\
   {\small \tt  weijinlong.hust@gmail.com }\\
   {\small \it $^d$ Department of Mathematics, Swansea University, Swansea SA2 8PP, UK}\\
   {\small \tt  j.l.wu@swansea.ac.uk }
}

\begin{document}
\maketitle

\medskip

\begin{abstract}
In this paper, we are concerned with regularity of nonlocal stochastic partial differential
equations of parabolic type. By using Companato estimates and Sobolev embedding theorem,
we first show the H\"{o}lder continuity (locally in the whole state space $\mathbb{R}^d$) for mild
solutions of stochastic nonlocal diffusion equations in the sense that the solutions $u$ belong to the
space $C^{\gamma}(D_T;L^p(\Omega))$ with the optimal H\"{o}lder continuity index $\gamma$
(which is given explicitly), where $D_T:=[0,T]\times D$ for $T>0$, and $D\subset\mathbb{R}^d$ being
a bounded domain. Then, by utilising tail estimates, we are able to obtain the estimates of mild
solutions in $L^p(\Omega;C^{\gamma^*}(D_T))$. What's more, we give an explicit
formula between the two index $\gamma$ and $\gamma^*$. Moreover, we prove H\"{o}lder continuity for mild
solutions on bounded domains. Finally, we present a new criteria to justify H\"{o}lder continuity for
the solutions on bounded domains. The novelty of this paper is that our method are suitable to the case
of time-space white noise.



\end{abstract}

\baselineskip=15pt

\section{Introduction}
\setcounter{equation}{0}

Given $T>0$ and $D\subset\mathbb{R}^d$, let $D_T:=[0,T]\times D$. Let $(\Omega,\mathcal{F},\{\mathcal{F}_t\}_{t\ge0},\mathbb{P})$ be a given
filtered probability space. In our previous paper \cite{LGWW}, we obtained regularity of singular stochastic
integrals in the following space
    \bess
\mathscr{L}^{p,\theta}((D_T;\delta);L^p(\Omega))
  \eess
for $p>1,\theta>0,\delta>0$. Further, by virtue of the celebrated Sobolev embedding theorem $\mathscr{L}^{p,\theta}(D;\delta)\hookrightarrow C^{\gamma}(\bar D;\delta)$
for $\theta>1$, we succeeded in obtaining estimates of solutions in the H\"older space
   \bess
C^{\gamma}(D_T;L^p(\Omega)),
  \eess
where $\gamma=\frac{(d+2)(\theta-1)}{p}$.
In the present paper, we aim to obtain the estimates of solutions in the space
   \bess
L^p(\Omega;C^{\gamma}(D_T)).
   \eess
The fundamental difficulty is the fact that usually
   \bess
\mathbb{E}\sup_{t,x}\neq
\sup_{t,x}\mathbb{E}.
   \eess
In this  paper, we are going to use the tail estimates to overcome the above mentioned difficulty.
The idea is fairly easy to explicate. In fact, note that
   \bess
\mathbb{E}(|X|^p)&=&\int_\Omega |X|^pd\mathbb{P}(\omega)\\
&=&p\int_0^\infty \mathbb{P}\{|X|>a\}a^{p-1}da\\
&=&p\int_0^M \mathbb{P}\{|X|>a\}a^{p-1}da+p\int_M^\infty \mathbb{P}\{|X|>a\}a^{p-1}da\\
&\leq&M^p+p\int_M^\infty \mathbb{P}\{|X|>a\}a^{p-1}da
   \eess
for any arbitrarily fixed constant $M>0$. In order to obtain the $L^p$-boundedness, by the above inequality,
we only need to show that the second integral is bounded. Further, by utilising Chebyshev's inequality, one
can derive the desired results by means of the estimates in $\mathscr{L}^{p,\theta}((D_T;\delta);L^p(\Omega))$.

Let us recall some regularity results about stochastic partial differential equations (SPDEs).
The earliest results about the $L_p$-theory of SPDEs appeared in the works of Krylov \cite{Krylov,K1996}.
Recently, Kim-Kim \cite{KK2016} considered the
$L_p$-theory for SPDEs
driven by L\'{e}vy processes, also see \cite{DMS,KKL,K2004,K2008}.
Zhang \cite{Z2006} obtained the $L_p$-theory of semi-linear SPDEs on general measure spaces.
 Let us also mention Zhang \cite{Z2013} where very interestingly $L_p$-maximal regularity of (deterministic) nonlocal parabolic PDEs and Krylov estimate for
SDEs driven by Cauchy processes are proved.

The H\"{o}lder estimate of SPDEs
has been studied by many authors. Let us mention a few. Hsu-Wang-Wang \cite{HW} established
the stochastic De Giorgi iteration and regularity of semilinear SPDEs. Du-Liu \cite{DL}
obtained the Schauder estimate for SPDEs. Combining the deterministic theory and convolution properties, Debussche-de Moor-Hofmanov\'a \cite{DMH}
established the regularity result for quasilinear SPDEs of parabolic type. Kuksin-Nadirashvili-Piatnitski \cite{KNP} obtained
H\"{o}lder estimates for solutions of parabolic SPDEs on bounded domains. Most recently,
Tian-Ding-Wei \cite{TDL} derived the local H\"{o}lder estimates of mild solutions of stochastic
nonlocal diffusion equations by using tail estimates \cite{KNP}. The results on H\"{o}lder estimate of PDEs with time-space white
noise are few. Fortunately, our method is suitable the  time-space white case.

There are two methods to deal with the Schauder estimate for SPDEs. One is using the smooth
property of kernel, the other is using the iteration technique. In this paper, we use the Morrey-Campanato
estimates and tail estimates to obtain the desired results.
The advantage of Morrey-Campanato estimates is to use the properties of kernel function and Sobolev
embedding theorem. Comparing with other methods to obtain the H\"{o}lder estimate, it is clear that this method
is relatively simple.

The rest of this paper is organized as follows. Section 2 presents some preliminaries. In section 3, we state and prove our main results
 on H\"older estimate over the whole spatial space. Section 4 is concerned with  H\"older estimate on bounded domains.
Section 5 is devoted to some applications of our main results.

\section{Preliminaries}
\setcounter{equation}{0}

Set, for $X=(t,x)\in\mathbb{R}\times\mathbb{R}^{d}$ and $Y=(s,y)\in\mathbb{R}\times\mathbb{R}^{d}$,
the following
   \bess
\delta(X,Y):=\max\left\{|x-y|,\,|t-s|^{\frac{1}{2}}\right\}.
   \eess
Let $Q_c(X)$ be the ball centered in $X=(t,x)$ with radius $c>0$, i.e.,
   \bess
Q_c(X):=\{Y=(s,y)\in\mathbb{R}\times\mathbb{R}^{d}:\,\delta(X,Y)<R\}
=(t-c^2,t+c^2)\times B_c(x).
  \eess
Fix $T\in(0,\infty)$ arbitrarily. Denote
   \bess
\mathcal {O}_T:=(0,T)\times\mathbb{R}^d.
  \eess
For  a bounded domain $D\subset\mathbb{R}^{d}$, we denote $D_T:=[0,T]\times D$. For a point $X\in D_T,
D(X,r):=D_T\cap Q_r(X)$ and $d(D):=diam (D)$ (that is, the diameter of $D$). Let us first give the definition of
Campanato space.

\begin{defi}\lbl{d6.1} (Campanato Space) Let $p\geq1$ and $\theta\geq0$.
The Campanato space $\mathscr{L}^{p,\theta}(D;\delta)$ is a subspace of
$L^p(D)$ such that
   \bess
[u]_{\mathscr{L}^{p,\theta}(D_T;\delta)}:=\left(\sup_{X\in D_T,d(D)\geq\rho>0}\frac{1}{|D(X,\rho)|^\theta}
\int_{D(X,\rho)}|u(Y)-u_{X,\rho}|^pdY\right)^{1/p}<\infty,\,\, u\in L^p(D_T)
    \eess
where $|D(X,\rho)|$ stands for the Lebesgue measure of the Borel set $D(X,\rho)$ and
   \bess
u_{X,\rho}:=\frac{1}{|D(X,\rho)|}\int_{D(X,\rho)}u(Y)dY.
   \eess
For $u\in \mathscr{L}^{p,\theta}(D_T;\delta)$, we define
    \bess
\|u\|_{\mathscr{L}^{p,\theta}(D_T;\delta)}
:=\left(\|u\|_{L^p(D_T)}^p+[u]_{\mathscr{L}^{p,\theta}(D_T;\delta)}^p\right)^{1/p}.
  \eess
\end{defi}

Next, we recall the definition of H\"{o}lder space.

\begin{defi}\lbl{d6.2} (H\"{o}lder Space) Let $0<\alpha\leq1$. A function $u$ belongs
to the H\"{o}lder space $C^{\alpha}(\bar D_T;\delta)$ if $u$ satisfies the following condition
   \bess
[u]_{C^{\alpha}(\bar D_T;\delta)}:=\sup_{X\in D_T,d(D)\geq\rho>0}\frac{|u(X)-u(Y)|}{\delta(X,Y)^\alpha}
<\infty.
    \eess
For $u\in C^{\alpha}(\bar D_T;\delta)$, we define
    \bess
\|u\|_{C^{\alpha}(\bar D_T;\delta)}:=\sup_{D_T}|u|+[u]_{C^{\alpha}(\bar D_T;\delta)}.
  \eess
\end{defi}

\begin{defi}\lbl{d6.3} Let $D_T\subset\mathbb{R}^{d+1}$ be a domain. We call the
domain $D_T$ an $A$-type domain if there exists a constant $A>0$ such that
$\forall X\in D_T$ and $\forall\, 0<\rho\leq d(D)$, it holds that
   \bess
|D_T(X,\rho)|=|D_T\cap Q_\rho(X)|\geq A|Q_\rho(X)|.
  \eess
\end{defi}
Recall that given two sets $B_1$ and $B_2$, the relation $B_1\cong B_2$ means
that both $B_1\subseteq B_2$ and $B_2\subseteq B_1$ hold. The notation $f(x)\approx g(x)$
means that there is a number $0<C<\infty$ independent of $x$,
i.e. a constant, such that for every $x$ we have $C^{-1}f(x)\leq g(x)\leq Cf(x)$. We have then the following
relation of the comparison of the two spaces defined above

\begin{prop}\lbl{p2.1} Assume that $D_T$ is an $A$-type bounded domain.
Then, for $p\geq1$ and for $1<\theta\leq1+\frac{p}{d+2}$ (Recall that $d$ is
the dimension of the space),
   \bess
\mathscr{L}^{p,\theta}(D_T;\delta)\cong C^{\gamma}(\bar D_T;\delta)
  \eess
with
   \bess
\gamma=\frac{(d+2)(\theta-1)}{p}.
   \eess
\end{prop}

We want to use the tail estimate to deive the following boundedness results
   \bess
\mathbb{E}\|u\|^p_{C^\gamma([0,T]\times D)}\leq C, \quad \forall p\geq 1
   \eess
for solutions $u$ of SPDEs. To this end, we need the following

  \begin{prop}\lbl{p2.2} \cite[Lemma 2.1]{TDL} Let $\rho_0\in L^p(\mathbb{R}^d\times\Omega)$.
Consider the Cauchy problem
   \bes
\partial_t\rho(t,x)=\Delta^\alpha\rho(t,x), \ \ t>0,\ x\in \mathbb{R}^d; \ \ \rho(0,x)=\rho_0(x).
  \lbl{2.1}\ees
Then, for any $0<\beta<1$, the following estimates for the unique mild solution of (\ref{2.1})
   \bes
\|\rho(t,\cdot)\|_{C^\beta(\mathbb{R}^d)}\leq Ct^{-\frac{\beta}{2\alpha}-\frac{d}{2p\alpha}}
\|\rho_0\|_{L^p(\mathbb{R}^d)}, \ \ \mathbb{P}- a.s. \,\,\omega\in\Omega,
   \lbl{2.2}\ees
and
    \bes
|\rho(t+\delta,x)-\rho(t,x)|\leq Ct^{-\beta-\frac{d}{2p\alpha}}\|\rho_0\|_{L^p(\mathbb{R}^d)},  \ \ \mathbb{P}- a.s. \,\,\omega\in\Omega.
   \lbl{2.3}\ees
  \end{prop}

We end this section with the following properties of kernel function $K$ satisfying $K_t=\Delta^\alpha K$ (the reader is referred to \cite{B-J, B-S-S, CH, I} for more details)
 \begin{itemize}
 \item for any $t>0$,
\bess
\|K(t, \cdot)\|_{L^1(\mathbb{R}^d)}=1 \text{ for all } t>0.
\eess
 \item $ K(t, x, y)$ is  $C^\infty$ on $(0,\infty)\times \mathbb{R}^d\times \mathbb{R}^d$ for each $t>0$;

 \item   for $t>0$, $x, y\in \mathbb{R}^d$, $x\neq y$, the sharp estimate of $K(t, x)$ is
\bess
K(t, x, y)\approx \min\left( \frac{t}{|x-y|^{d+2\alpha}}, t^{-d/(2\alpha)}\right);
  \eess
\item  for $t>0$, $x, y\in \mathbb{R}^d$, $x\neq y$,   the  estimate of the first order derivative of $ K(t, x)$ is
\bes
  |\nabla_x K(t, x, y)|\approx  |y-x|\min\left\{ \frac{t}{|y-x|^{d+2+2\alpha}}, t^{-\frac{d+2}{2\alpha}}\right\}.
 \lbl{4.4}\ees
\end{itemize}

The estimate (\ref{4.4}) for the first order derivative of $K(t,x)$ was derived in  \cite[Lemma 5]{B-J}.
Xie et al. \cite{XDLL} obtained the estimate of the $m$-th order derivative of $p(t, x)$ by induction.
\begin{prop}
\label{p2.3}{\rm\cite[Lemma 2.1]{XDLL}}
For any $m\geq 0$, we have
  \bess
\partial_x^m K(t, x)=\sum_{n=0}^{n=\lfloor \frac{m}{2}\rfloor}C_n |x|^{m-2n} \min \left\{ \frac{t}{|x|^{d+2\alpha+2(m-n)}}, t^{-\frac{d+2(m-n)}{2\alpha}}\right\},
    \eess
where $\lfloor \frac m 2\rfloor$ means the largest integer that is less  than $\frac m2$.
\end{prop}

\section{H\"{o}lder estimate locally over the whole spatial space}
\setcounter{equation}{0}
In this section, we establish the Morrey-Campanato estimates
under different assumption on stochastic term.
Set
 \bess
\mathcal {K}g(t,x):=\int_0^t\int_{\mathbb{R}^d}K(t-r,y)g(r,x-y)dydW(r).
  \eess

The first result is similar to the deterministic case.
We consider the following equation
   \bes
du_t=\Delta^\alpha udt+g(t,x)dW_t,\ \ \ u|_{t=0}=0,
   \lbl{3.1}\ees
where $\Delta^\alpha=-(-\Delta)^\alpha$ and $W_t$ is a
standard Brownian motion on a filtered probability space
$(\Omega,\mathcal {F},\{\mathcal {F}_t\}_{t\ge0},\mathbb{P})$.
\begin{theo}\lbl{t3.1} Let $D$ be an $A$-type bounded domain in $\mathbb{R}^{d+1}$ such that $\bar D\subset \mathcal {O}_T$.
Suppose that $g\in L^\infty_{loc}(\mathbb{R}_+;L^p(\Omega\times\mathbb{R}^d))$ for $p>d/\alpha$
 is $\mathcal {F}_t$-adapted process, and that $0<\beta<\alpha$ satisfies $(\alpha-\beta) p-d\geq0$. Then, there is a mild solution $u$ of (\ref{3.1}) and
$u\in \mathscr{L}^{p,\theta}((D_T;\delta);L^p(\Omega))\cap L^p(\Omega;C^{\beta}(D_T))$. Moreover, it holds that
   \bes
&&\|u\|_{\mathscr{L}^{p,\theta}((D_T;\delta);L^p(\Omega))}\leq C\|g\|_{ L^\infty([0,T];L^p(\Omega\times\mathbb{R}^d))},\lbl{3.2}\\
&&\|u\|_{C^{\beta}(D_T;L^p(\Omega))}\leq C\|g\|_{ L^\infty([0,T];L^p(\Omega\times\mathbb{R}^d))},\lbl{3.3}
  \ees
where $\theta=1+\frac{\beta p}{d+2}$. Moreover,  taking $0<\delta<\beta p/2$ and $q>(d+2)/\delta$,
we have for $0<r<q$
   \bes
\|u\|_{L^r(\Omega;C^{\beta^*}(D_T))}\leq C\|g\|_{ L^\infty([0,T];L^p(\Omega\times\mathbb{R}^d))},\lbl{3.4}
   \ees
where $\beta^*=\beta-2\delta/p$.
  \end{theo}

{\bf Proof.}  The existence of mild solution of (\ref{3.1}) is a classical result under
the above assumptions. Now we prove the inequality (\ref{3.2}).
Due to the definition of Companato space, it suffices to show that
   \bess
[u]_{\mathscr{L}^{p,\theta}((D_T;\delta);L^p(\Omega))}<\infty.
  \eess
Direct calculus shows that
    \bess
[u]^p_{\mathscr{L}^{p,\theta}((D_T;\delta);L^p(\Omega))}&\leq&
\sup_{D(X,c),X\in D_T,0<c\leq d(D)}\frac{1}{|D(X,c)|^{1+\theta}}\\
&&\times\mathbb{E}
\int_{D(X,c)}\int_{D(X,c)}|u(t,x)-u(s,y)|^pdtdxdsdy\\
&\leq&\sup_{D(X,c),X\in D_T,0<c\leq d(D)}\frac{1}{|D(X,c)|^{1+\theta}}\\
&&\times\mathbb{E}
\int_{D(X,c)}\int_{D(X,c)}\Big|\int_0^t\int_{\mathbb{R}^d}K(t-r,x-z)g(r,z)dzdW(r)\\
&&-\int_0^s\int_{\mathbb{R}^d}K(s-r,y-z)g(r,z)dzdW(r)\Big|^p\\
&:=&\sup_{D(X,c),X\in D_T,0<c\leq d(D)}\frac{1}{|D(X,c)|^{1+\theta}}
\int_{D(X,c)}\int_{D(X,c)}\mathbb{E}\Upsilon dtdxdsdy.
   \eess

Set $t\geq s$. We have the following estimates
    \bess
\mathbb{E}\Upsilon
&\leq&C\mathbb{E}\Big|\int_0^s\int_{\mathbb{R}^d}(K(t-r,x-z)-K(s-r,y-z))g(r,z)dzdW(r)\Big|^p\\
&&
+C\mathbb{E}\Big|\int_s^t\int_{\mathbb{R}^d}K(t-r,x-z)g(r,z)dzdW(r)\Big|^p\\
&\leq&C\mathbb{E}\Big|\int_0^s\left(\int_{\mathbb{R}^d}(K(t-r,x-z)-K(s-r,y-z))g(r,z)dz\right)^2dr\Big|^{\frac{p}{2}}\\
&&
+C\mathbb{E}\Big|\int_s^t\left(\int_{\mathbb{R}^d}K(t-r,x-z)g(r,z)dz\right)^2dr\Big|^{\frac{p}{2}}\\
&:=&C(H_1+H_2).
   \eess

{\bf Estimate of $H_1$}.

Take $\beta>0$ satisfying $(\alpha-\beta)p-d\geq0$.
We first recall the following fractional mean value formula
(see (4.4) of \cite{Jumarie2006})
   \bess
f(x+h)=f(x)+\Gamma^{-1}(1+\beta)h^\beta f^{(\beta)}(x+\theta h),
   \eess
where $0<\beta<1$ and $0\leq\theta\leq1$ depends on $h$ satisfying
   \bess
\lim\limits_{h\downarrow0}\theta^\beta=\frac{\Gamma^2(1+\beta)}{\Gamma(1+2\beta)},
   \eess
By using the Propositions \ref{p2.2} and \ref{p2.3}, the above fractional mean value formula and H\"{o}lder inequality, we have
   \bess
H_1 &=&\mathbb{E}\Big|\int_0^s\left(\int_{\mathbb{R}^d}(K(t-r,x-z)-K(s-r,y-z))g(r,z)dz\right)^2dr\Big|^{\frac{p}{2}}\\
 &\leq&C\mathbb{E}\Big|\int_0^s\left(\int_{\mathbb{R}^d}|K(t-r,x-z)-K(s-r,x-z)|\cdot|g(r,z)|dz\right)^2dr\Big|^{\frac{p}{2}}\\
 &&+C\mathbb{E}\Big|\int_0^s\left(\int_{\mathbb{R}^d}(K(s-r,x-z)-K(s-r,y-z))\cdot g(r,z) dz\right)^2dr\Big|^{\frac{p}{2}}\\
 &\leq&C(t-s)^{\frac{\beta p}{2}}\mathbb{E}\Big|\int_0^s\left(\int_{\mathbb{R}^d}
 |\frac{\partial^{\frac{\beta }{2}} K}{\partial t^{\frac{\beta }{2}}}(\xi-r,x-z)|^qdz \right)^{\frac{2}{q}}
\|g(r)\|^2_{L^p(\mathbb{R}^d)}dr\Big|^{\frac{p}{2}}\\
&&+C|x-y|^{\beta p}\mathbb{E}\Big|\int_0^s (s-r)^{-\frac{\beta}{\alpha}-\frac{d}{p\alpha}}
\|g(r)\|^2_{L^p(\mathbb{R}^d)}dr\Big|^{\frac{\delta}{2}}\\
&\leq&C(t-s)^{\frac{\beta p}{2}}\|g\|_{L^p(\Omega;L^\infty([0,T];L^p(\mathbb{R}^d)))}^p
\left[\int_0^s\left(\int_{\mathbb{R}^d}|\frac{\partial^{\frac{\beta }{2}} K}{\partial t^{\frac{\beta }{2}}}(\xi-r,x-z)|^qdz \right)^{\frac{2}{q}}dr\right]^{\frac{p}{2}}\\
&&+C|x-y|^{\beta p}\|g\|_{L^p(\Omega;L^\infty([0,T];L^p(\mathbb{R}^d)))}^p\left[\int_0^s(s-r)^{-\frac{\beta}{\alpha}-\frac{d}{p\alpha}}dr\right]^{\frac{p}{2}}\\
&\leq&C((t-s)^{\frac{\beta p}{2}}+|x-y|^{\beta p}),
   \eess
where $q=p/(p-1)$, $\xi=\theta t+(1-\theta)s$, and we used the
following fact
   \bess
&&\int_0^s\left(\int_{\mathbb{R}^d}|\frac{\partial^{\frac{\beta }{2}} K}{\partial t^{\frac{\beta }{2}}}(\xi-r,x-z)|^qdz \right)^{\frac{2}{q}}dr
\\
&\leq&C\int_0^s\left(\int_0^{(\xi-r)^{\frac{1}{2\alpha}}}(\xi-r)^{-\frac{dq+2q\alpha\beta}{2\alpha}}|z|^{q\alpha\beta+d-1}d|z|\right.\\
&&\left.
+\int_{(\xi-r)^{\frac{1}{2\alpha}}}^\infty(\xi-r)^q|z|^{-(qd+2q\alpha+2q\alpha\beta)}|z|^{q\alpha\beta+d-1}d|z| \right)^{\frac{2}{q}}dr\\
&\leq&C\left[(\theta(t-s))^{\frac{d-dq+q\alpha(1-\beta)}{q\alpha}}+\xi^{\frac{d-dq+q\alpha(1-\beta)}{q\alpha}}\right]\\
&\leq&C
  \eess
because using $q=p/(p-1)$, we have
   \bess
d-dq+q\alpha(1-\beta)>0\Leftrightarrow p(\alpha-\alpha\beta)>d\Leftarrow p(\alpha-\beta)>d.
  \eess

Similarly, we have
  \bess
\int_0^s(s-r)^{-\frac{\beta}{\alpha}-\frac{d}{p\alpha}}dr=\frac{p\alpha}{(\alpha-\beta)p-d}s^{\frac{(\alpha-\beta)p-d}{p\alpha}}
\leq C
   \eess
provided that $(\alpha-\beta)p-d\geq0$.

{\bf Estimate of $H_2$}.

Similar to the estimate of $H_1$, we have
   \bess
H_2&=&\mathbb{E}\Big|\int_s^t\left(\int_{\mathbb{R}^d}K(t-r,x-z)g(r,z)dz\right)^2dr\Big|^{\frac{p}{2}}\\
&\leq&\|g\|_{L^p(\Omega;L^\infty([0,T];L^p(\mathbb{R}^d)))}^p
\left[\int_s^t\left(\int_{\mathbb{R}^d}|K(t-r,x-z)|^qdz \right)^{\frac{2}{q}}dr\right]^{\frac{p}{2}}\\
&\leq&C\|g\|_{L^p(\Omega;L^\infty([0,T];L^p(\mathbb{R}^d)))}^p(t-s)^{\frac{q\alpha-(q-1)d}{q\alpha}\times \frac{p}{2}}
  \eess
provided that $\alpha p> d$. Indeed, by using $1/p+1/q=1$, we have
   \bess
q\alpha-(q-1)d>0\ \ \Longleftrightarrow \ \ \alpha p> d.
  \eess

Combining the assumption of $p$, we have
    \bess
H_2\leq C(t-s)^{\frac{p\alpha-d}{2\alpha}}.
   \eess
Assume that $D(X,c)=D_T\cap Q_c$ and $Q_c=Q_{c}(t_0,x_0)$.
Noting that $(t,x)\in Q_{c}(t_0,x_0)$ and $(s,y)\in Q_{c}(t_0,x_0)$, we have
   \bess
0\leq t-s\leq 2{c}^2\ \ \ {\rm and}\ \ \ |x-y|\leq|x-x_0|+|y-x_0|\leq2{c}.
   \eess
By using the definition of $A$-type bounded domain, we have
       \bess
[u]_{\mathscr{L}^{p,\theta}((D_T;\delta);L^p(\Omega))}&\leq&
\sup_{D(X,c),X\in D_T,0<c\leq d(D)}\frac{1}{|D(X,c)|^{1+\theta}}\mathbb{E}
\int_{D(X,c)}\int_{D(X,c)}\mathbb{E}\Upsilon dtdxdsdy\\
&\leq&C\|g\|_{L^p(\Omega;L^\infty([0,T];L^p(\mathbb{R}^d)))}^p,
   \eess
where $\theta=1+\frac{\beta p}{d+2}$. This yields the inequality (\ref{3.2}).
Applying Proposition \ref{p2.1}, one can obtain the inequality (\ref{3.3}).

Next, we prove the inequality (\ref{3.4}).
In order to use the technique of tail estimates, we first consider the
following estimates.
Let $(t_0,x_0)\in D_T\subset\mathcal {O}_T$ and
  \bess
Q_c(t_0,x_0)=(t_0-c^2,t_0+c^2)\times B_c(x_0).
   \eess
Then we have $\bar D_T\subset Q_{d(D)}(t_0,x_0)$.
Set $(t_1,x_1),(t_2,x_2)\in D_T$ , $Q_i:=D_T\cap Q_{c_i}(t_i,x_i)$, $i=1,2$
and
   \bess
F(t_i,x_i,c_i)&=& \frac{1}{|Q_i|^{1+\theta}}
\int_{Q_i}\int_{Q_i}|u(t,x)-u(s,y)|^pdtdxdsdy\\
&=&\frac{1}{|Q_i|^{1+\theta}}
\int_{Q_i}\int_{Q_i}|\mathcal {K}g(t,x)-\mathcal {K}g(s,y)|^pdtdxdsdy.
   \eess
Notice that
   \bess
F(t_1,x_1,c_1)-F(t_2,x_2,c_2)&=&[F(t_1,x_1,c_1)-F(t_2,x_1,c_1)]\\
&&+[F(t_2,x_1,c_1)-F(t_2,x_2,c_1)]\\
&&+[F(t_2,x_2,c_1)-F(t_2,x_2,c_2)]\\
&:=&I_1+I_2+I_3.
  \eess

Estimate of $I_1$:
   \bess
I_1&=&F(t_1,x_1,c_1)-F(t_2,x_1,c_1)\\
&=&\frac{1}{|Q_1|^{1+\theta}}
\int_{Q_1}\int_{Q_1}|\mathcal {K}g(t,x)-\mathcal {K}g(s,y)|^pdtdxdsdy\\
&&
-\frac{1}{|Q_{12}|^{1+\theta}}
\int_{Q_{12}}\int_{Q_{12}}|\mathcal {K}g(t,x)-\mathcal {K}g(s,y)|^pdtdxdsdy\\
&=&\frac{1}{|Q_1|^{1+\theta}}\left\{
\int_{Q_1\setminus Q_{12}}\int_{Q_1\setminus Q_{12}}|\mathcal {K}g(t,x)-\mathcal {K}g(s,y)|^pdtdxdsdy\right.\\
&&\left.+\int_{Q_{12}\setminus Q_{1}}\int_{Q_{12}\setminus Q_{1}}|\mathcal {K}g(t,x)-\mathcal {K}g(s,y)|^pdtdxdsdy\right\}\\
&&
+\left[\frac{1}{|Q_{1}|^{1+\theta}}
-\frac{1}{|Q_{12}|^{1+\theta}}\right]
\int_{Q_{12}}\int_{Q_{12}}|\mathcal {K}g(t,x)-\mathcal {K}g(s,y)|^pdtdxdsdy\\
&:=&I_{11}+I_{12},
   \eess
where $Q_{12}=D_T\cap Q_{c_1}(t_2,x_1)$.
For simplicity, we assume that $|Q_1|\geq|Q_{12}|$. Otherwise, we can chance the place of
$Q_1$ and $Q_{12}$. And thus $I_{12}\leq0$ almost surely. Now, we consider the term
$I_{11}$. Before giving the estimates of $I_{11}$, we first recall our aim. In order to
apply the tail estimate, we want to obtain
the estimates of $I_{11}$ like the followings:
   \bess
\mathbb{E}I_{11}\leq C(t_1-t_2)^\delta\ \ \ {\rm for}\ \ {\rm some}\ \ \delta>0.
   \eess
It is easy to see that
   \bess
|Q_1\setminus Q_{12}|\leq C(t_1-t_2)c_1^{d}\ \ {\rm and}\ \
|Q_1|\approx Cc_1^{d+2}.
  \eess
So we must put some assumption on $g$ in order to get some help from it.

Set $t>s$. Denote
  \bess
&&\mathbb{E}\int_{Q_1\setminus Q_{12}}\int_{Q_1\setminus Q_{12}}|\mathcal {K}g(t,x)-\mathcal {K}g(s,y)|^pdtdxdsdy\\
&=&\mathbb{E}\int_{Q_1\setminus Q_{12}}\int_{Q_1\setminus Q_{12}}\mathbb{E}\Upsilon dtdxdsdy.
  \eess

Similar to the proof of inequality (\ref{3.2}), we have
   \bess
\mathbb{E}\Upsilon\leq Cc_1^{\beta p}.
   \eess
Noting that $(t,x)\in Q_{1}$ and $(s,y)\in Q_{1}$, we have
   \bess
0\leq t-s\leq 2{c_1}^2\ \ \ {\rm and}\ \ \ |x-y|\leq|x-x_1|+|y-x_1|\leq2{c_1}.
   \eess
Using the above inequalities and the properties of $A$-type domain, we deduce
  \bess
&&\mathbb{E}\int_{Q_1\setminus Q_{12}}\int_{Q_1\setminus Q_{12}}\mathbb{E}\Upsilon dtdxdsdy\\
&\leq& C(p,T)c_1^{\beta p}|Q_1\setminus Q_{12}|^2\|g\|_{L^p(\Omega;L^\infty([0,T];L^p(\mathbb{R}^d)))}^p.
   \eess
Since $D_T$ is a $A$-type bounded domain, we have for $2c_1\leq diam D$,
   \bess
&&A|Q_{c_1}(t_1,x_1)|\leq |Q_1|\leq |Q_{c_1}(t_1,x_1)|\\
&&A|Q_{c_1}(t_1,x_1)\setminus Q_{c_1}(t_2,x_1)|\leq |Q_1\setminus Q_{12}|\leq |Q_{c_1}(t_1,x_1)\setminus Q_{c_1}(t_2,x_1)|.
   \eess
We remark that
   \bess
&&|Q_{c_1}(t_1,x_1)|\approx Cc_1^{d+2},\\
&&|Q_{c_1}(t_1,x_1)\setminus Q_{c_1}(t_2,x_1)|\leq Cc_1^{d}[c_1^2\wedge(t_1-t_2)],
   \eess
where $C$ is a positive constant which does not depend on
$c_1$. Noting that $Q_1\setminus Q_{12}\subset Q_1$ and taking $0<\delta<\beta p/2$, we have
  \bess
&&\mathbb{E}
\int_{Q_1\setminus Q_{12}}\int_{Q_1\setminus Q_{12}}|\mathcal {K}g(t,x)-\mathcal {K}g(s,y)|^pdtdxdsdy\\
&\leq&C(C_0,D,d,T)\|g\|_{L^p(\Omega;L^\infty([0,T];L^p(\mathbb{R}^d)))}^p|Q_1|^{2+\frac{\beta p-2\delta}{d+2}}|t_1-t_2|^\delta.
   \eess
Similarly, we can get
 \bess
&&\mathbb{E}\int_{Q_{12}\setminus Q_{1}}\int_{Q_{12}\setminus Q_{1}}|u(t,x)-u(s,y)|^pdtdxdsdy\\
&\leq&C(D,d,T)\|g\|_{L^p(\Omega;L^\infty([0,T];L^p(\mathbb{R}^d)))}^p|Q_1|^{2+\frac{\beta p-2\delta}{d+2}}|t_1-t_2|^\delta.
   \eess
Due to the fact that $I_{12}\leq0$, we have
   \bess
\mathbb{E}I_1\leq C(D,d,T)\|g\|_{L^p(\Omega;L^\infty([0,T];L^p(\mathbb{R}^d)))}^p|t_1-t_2|^\delta,
  \eess
where $\theta=1+\frac{\beta p-2\delta}{d+2}$.

Next, we estimate $I_2$. By using the fact that
   \bess
\Big|[D\cap Q_{c_1}(t_2,x_1)]\setminus [D\cap Q_{c_1}(t_2,x_2)]\Big|\leq Cc_1^{d+1}|x_1-x_2|,
  \eess
similar to the estimates of $I_1$, we can take $0<\delta<\beta p/2$ such that
    \bess
\mathbb{E}I_2&=&\mathbb{E}[F(t_2,x_1,c_1)-F(t_2,x_2,c_1)]\\
&\leq& C(D,d,T)\|g\|_{L^p(\Omega;L^\infty([0,T];L^p(\mathbb{R}^d)))}^p|x_1-x_2|^\delta,
   \eess
where $\theta=1+\frac{\beta p-\delta}{d+2}$.

Next, we estimate $I_3$. By using the fact that
   \bess
\Big|[D\cap Q_{c_1}(t_2,x_2)]\setminus [D\cap Q_{c_2}(t_2,x_2)]\Big|\leq Cc_1^{d+1}(c_1-c_2), \ \ {\rm if} \ \ c_1\geq c_2,
  \eess
similar to the estimates of $I_1$, we can estimate
    \bess
\mathbb{E}I_3&=&\mathbb{E}[F(t_2,x_2,c_1)-F(t_2,x_2,c_2)]\\
&\leq& C(D,d,T)\|g\|_{L^p(\Omega;L^\infty([0,T];L^p(\mathbb{R}^d)))}^p|c_1-c_2|^\delta,
   \eess
where $\theta=1+\frac{\beta p-\delta}{d+2}$.

Therefore, we have
   \bess
 &&\mathbb{E}|F(t_1,x_1,c_1)-F(t_2,x_2,c_2)|^q\\
 &\leq& C(D,d,T)\|g\|_{L^p(\Omega;L^\infty([0,T];L^p(\mathbb{R}^d)))}^{pq}(|t_1-t_2|+|x_1-x_2|+|c_1-c_2|)^{\delta q},
    \eess
where $\theta=1+\frac{\beta p-2\delta}{d+2}$, $(t_i,x_i)\in D_T$ and $0<c_i\leq d(D)$, $i=1,2$.

For simplicity, we set $D_T=[0,1]^{d+1}$ and $c\in[0,2]$. One introduces a
sequence of sets:
   \bess
\mathcal{S}_n=\{z\in\mathbb{Z}^{d+2}| \ z2^{-n}\in(0,1)^{d+1}\times(0,2)\}, \ \ n\in \mathbb{N}.
   \eess
For an arbitrary $e=(e_1,e_2,\cdots,e_{d+2})\in\mathbb{Z}^{d+2}$ such that
   \bess
|e|_\infty=\max_{1\leq j\leq d+2}|e_j|=1,
   \eess
and for every $z,z+e\in \mathcal{S}_n$, we define $v_z^{n,e}=|F((z+e)2^{-n})-F(z2^{-n})|$.
From the above discussion, we have
   \bess
\mathbb{E}|v_z^{n,e}|^q\leq C(\beta,C_0,D,d,T)\|g\|_{L^p(\Omega;L^\infty([0,T];L^p(\mathbb{R}^d)))}^{pq}
2^{-n\delta q}:=\hat C2^{-n\delta q}.
   \eess
For any $\tau>0$ and $K>0$, one sets a number of events
   \bess
\mathcal{A}_{z,\tau}^{n,e}=\{\omega\in\Omega |  v_z^{n,e}\geq K\tau^n, \ z,z+e\in\mathcal{S}_n\},
  \eess
which yields that
    \bess
\mathbb{P}(\mathcal{A}_{z,\tau}^{n,e})\leq\frac{\mathbb{E}|v_z^{n,e}|^q}{K^q\tau^{qn}}
\leq\frac{\hat C2^{-n\delta q}}{K^q\tau^{qn}}.
  \eess
Noting that for each $n$, the total number of the events $\mathcal{A}_{z,\tau}^{n,e},z,z+e\in\mathcal{S}_n$
is not larger than $2^{d+2}3^{d+2}$. Hence the probability of the union
   \bess
\mathcal{A}_\tau^n=\cup_{z,z+e\in\mathcal{S}_n}(\cup_{\|e\|_\infty=1}\mathcal{A}_{z,\tau}^{n,e})
    \eess
meets the estimate
   \bess
\mathbb{P}(\mathcal{A}_\tau^n)
\leq\frac{\hat C2^{-n\delta q}}{K^q\tau^{qn}}2^{(d+2)n}\leq\hat CK^{-q}\left(\frac{2^{d+2}}{(2^\delta\tau)^q}\right)^n.
   \eess
Let $\tau=2^{-\nu\delta}$, where $\nu>0$ satisfies
$(1-\nu)\delta q\geq d+2$. Then the probability of the event $\mathcal{A}=\cup_{n\geq1}\mathcal{A}_\tau^n$
can be calculated that
   \bes
\mathbb{P}(\mathcal{A})
\leq C\hat CK^{-q}.
   \lbl{3.5}\ees
For every point $\xi=(t,x,c)\in(0,1)^{d+1}\times(0,2)$, we have
$\xi=\sum_{i=0}^\infty e_i2^{-i}$ ($\|e_i\|_\infty\leq1$). Denote
$\xi_k=\sum_{i=0}^ke_i2^{-i}$ and $\xi_0=0$. For any $\omega\notin\mathcal{A}$, we have
$|F(\xi_{k+1})-F(\xi_k)|<K\tau^{k+1}$, which implies that
   \bes
|F(t,x,c)|\leq\sum_{k=0}^\infty|F(\xi_{k+1})-F(\xi_k)|<K\sum_{k=1}^\infty\tau^k\leq K(2^{\nu\delta}-1)^{-1}.
   \lbl{3.6} \ees
Set $v_1=\sup_{(t,x,c)\in(0,1)^{d+1}\times(0,2)}|F(t,x,c)|$, then
$v_1=\sup_{(t,x,c)\in[0,1]^{d+1}\times[0,2]}|F(t,x,c)|$ since $F$ has a continuous
version. For $0<r<q$, we have
   \bes
\mathbb{E}v_1^r=r\int_0^\infty a^{r-1}\mathbb{P}(v_1\geq a)da
=r\int_0^{\gamma K}a^{r-1}\mathbb{P}(v_1\geq a)da+r\int_{\gamma K}^\infty a^{r-1}\mathbb{P}(v_1\geq a)da.
   \lbl{3.7}\ees
If one chooses $\gamma\geq (2^{\nu\delta}-1)^{-1}$, using (\ref{3.5}), (\ref{3.6}) and (\ref{3.7}), we get
    \bess
\mathbb{E}v_1^r&\leq& (\gamma K)^r+C\hat C^{q}r\int_{\gamma K}^\infty a^{r-1-q}da\\
&\leq& (\gamma K)^r+C\hat Cr(cK)^{r-q},
   \eess
which yields that
   \bess
\mathbb{E}v_1^r\leq C(D,d,T)\|g\|_{L^p(\Omega;L^\infty([0,T];L^p(\mathbb{R}^d)))}^{pr}
2^{-n\delta q},
   \eess
if we choose $K=\|g\|_{L^p(\Omega;L^\infty([0,T];L^p(\mathbb{R}^d)))}^{p}$.  By using the following
embed inequality
   \bess
L^p(\Omega;\mathscr{L}^{p,\theta}(D_T;\delta))
  \cong L^p(\Omega;C^{\gamma}(\bar D_T;\delta)),
  \eess
we obtain the inequality (\ref{3.4}).
 The proof is complete. $\Box$

\begin{remark}\lbl{r3.1} It follows from Theorem \ref{t3.1} that the index $\beta$ and
$\beta^*$ satisfy $\beta>\beta^*$, which implies that if we want to change the places of
$\mathbb{E}$ and $\sup_{t,x}$, we must pay it on the index.

Comparing with the earlier results of \cite{TDL} (Tian et al. obtained the H\"{o}lder estimate to equation
(\ref{3.1}) locally in $\mathbb{R}^d$), we find the H\"{o}lder continuous index in this paper is larger than
that in  \cite{TDL}. More precisely, we obtain the index of time variable is closed to $1/2$. Since the
index of H\"{o}lder continuous of Brownian motion is $\frac{1}{2}-$, maybe the index obtained in
this paper is optimal.
\end{remark}

Next, we consider another case. If $g$ is a H\"{o}lder continuous function, the following
theorem shows that what assumptions should be put on the kernel function $K$.

 \begin{theo}\lbl{t3.2}
Let $u=\mathcal{K}\ast g$ and $D_T$ be an $A$-type bounded domain in $\mathbb{R}^{d+1}$ such that
$\bar D_T\subset \mathcal {O}_T$.
Suppose that $g\in C^\beta(\mathbb{R}_+\times\mathbb{R}^d)$, $0<\beta<1$, is a non-random function
  and $g(0,0)=0$.
Assume that there exists positive constants $\gamma_i$ ($i=1,2$) such that
the non-random kernel function satisfies that for any $t\in(0,T]$
  \bes
&&\int_0^s\left(\int_{\mathbb{R}^d}|K(t-r,z)-K(s-r,z)|(1+|z|^\beta)dz\right)^2dr\leq C(T,\beta)(t-s)^{\gamma_1},\lbl{3.8}\\
&&\int_0^s\left(\int_{\mathbb{R}^d}|K(s-r,z)|dz\right)^2dr\leq C_0,\lbl{3.9}\\
&& \int_s^t\left(\int_{\mathbb{R}^d}|K(t-r,z)|(1+|z|^\beta) dz\right)^2dr\leq C(T,\beta)(t-s)^{\gamma_2},
   \lbl{3.10} \ees
where $C_0$ is a positive constant. Then we have, for
$p\geq1$ and $\beta<\gamma$,
    \bes
&&\|u\|_{\mathscr{L}^{p,\theta}((D_T;\delta);L^p(\Omega))}\leq
C\|g\|_{ C^\beta(\mathbb{R}_+\times\mathbb{R}^d))},\nonumber\\
&&\|u\|_{C^{\beta}(D_T;L^p(\Omega))}\leq C\|g\|_{ C^\beta(\mathbb{R}_+\times\mathbb{R}^d))},
  \lbl{3.11}\ees
where $\theta=1+\frac{\gamma p}{d+2}$ and $\gamma=\min\{\gamma_1,\gamma_2,\beta\}$.
Moreover,  taking $0<\delta<\gamma p/2$ and $q>(d+2)/\delta$,
we have for $0<r<q$
   \bes
\|u\|_{L^r(\Omega;C^{\beta^*}(D_T))}\leq C\|g\|_{ C^\beta(\mathbb{R}_+\times\mathbb{R}^d))},
   \lbl{3.12}\ees
where $\beta^*=\gamma-2\delta/p$.
\end{theo}

{\bf Proof.} The proof of the (\ref{3.11}) is contained in our paper \cite{LGWW}. And we only
focus on the proof of (\ref{3.12}).

Similar to the proof of Theorem \ref{t3.1}, we need to estimate $I_i$, $i=1,2,3$.
Estimate of $I_1$:
   \bess
I_1&=&F(t_1,x_1,c_1)-F(t_2,x_1,c_1)\\
&=&\frac{1}{|Q_1|^{1+\theta}}
\int_{Q_1}\int_{Q_1}|\mathcal {K}g(t,x)-\mathcal {K}g(s,y)|^pdtdxdsdy\\
&&
-\frac{1}{|Q_{12}|^{1+\theta}}
\int_{Q_{12}}\int_{Q_{12}}|\mathcal {K}g(t,x)-\mathcal {K}g(s,y)|^pdtdxdsdy\\
&=&\frac{1}{|Q_1|^{1+\theta}}\left\{
\int_{Q_1\setminus Q_{12}}\int_{Q_1\setminus Q_{12}}|\mathcal {K}g(t,x)-\mathcal {K}g(s,y)|^pdtdxdsdy\right.\\
&&\left.+\int_{Q_{12}\setminus Q_{1}}\int_{Q_{12}\setminus Q_{1}}|\mathcal {K}g(t,x)-\mathcal {K}g(s,y)|^pdtdxdsdy\right\}\\
&&
+\left[\frac{1}{|Q_{1}|^{1+\theta}}
-\frac{1}{|Q_{12}|^{1+\theta}}\right]
\int_{Q_{12}}\int_{Q_{12}}|\mathcal {K}g(t,x)-\mathcal {K}g(s,y)|^pdtdxdsdy\\
&:=&I_{11}+I_{12},
   \eess
where $Q_{12}=D\cap Q_{c_1}(t_2,x_1)$.
For simplicity, we assume that $|Q_1|\geq|Q_{12}|$. Otherwise, we can chance the place of
$Q_1$ and $Q_{12}$. And thus $I_{12}\leq0$ almost surely.

It is easy to see that
   \bess
|Q_1\setminus Q_{12}|\leq C(t_1-t_2)c_1^{d}\ \ {\rm and}\ \
|Q_1|\approx Cc_1^{d+2}.
  \eess
So we must put some assumption on $g$ in order to get some help from it.

Set $t>s$. By the BDG inequality, we have
  \bess
&&\mathbb{E}\int_{Q_1\setminus Q_{12}}\int_{Q_1\setminus Q_{12}}|\mathcal {K}g(t,x)-\mathcal {K}g(s,y)|^pdtdxdsdy\\
&=&\mathbb{E}\int_{Q_1\setminus Q_{12}}\int_{Q_1\setminus Q_{12}}\Big|\int_0^t\int_{\mathbb{R}^d}K(t-r,z)g(r,x-z)dzdW(r)\\
&&
-\int_0^s\int_{\mathbb{R}^d}K(s-r,z)g(r,y-z)dzdW(r)\Big|^pdtdxdsdy\\
&\leq&2^{p-1}\mathbb{E}\int_{Q_1\setminus Q_{12}}\int_{Q_1\setminus Q_{12}}\Big|\int_0^s\int_{\mathbb{R}^d}(K(t-r,z)-K(s-r,z))
g(r,x-z)dzdW(r)\Big|^p \\
&&+2^{p-1}\mathbb{E}\int_{Q_1\setminus Q_{12}}\int_{Q_1\setminus Q_{12}}\Big|\int_0^s\int_{\mathbb{R}^d}K(s-r,z)
(g(r,x-z)-g(r,y-z))dzdW(r)\Big|^p \\
&&+2^{p-1}\mathbb{E}\int_{Q_1\setminus Q_{12}}\int_{Q_1\setminus Q_{12}}\Big|\int_s^t\int_{\mathbb{R}^d}K(t-r,z)
g(r,x-z)dzdW(r)\Big|^pdtdxdsdy\\
&\leq&C(p)\int_{Q_1\setminus Q_{12}}\int_{Q_1\setminus Q_{12}}
\left(\int_0^s|\int_{\mathbb{R}^d}|K(t-r,z)-K(s-r,z)||g(r,x-z)|dz|^2
dr\right)^{\frac{p}{2}} \\
&&
+C(p)\int_{Q_1\setminus Q_{12}}\int_{Q_1\setminus Q_{12}}
\left(\int_0^s|\int_{\mathbb{R}^d}|K(s-r,z)||g(r,x-z)-g(r,y-z)|dz|^2
dr\right)^{\frac{p}{2}} \\
&&
+C(p)\int_{Q_1\setminus Q_{12}}\int_{Q_1\setminus Q_{12}}\left(\int_s^t|
\int_{\mathbb{R}^d}K(t-r,z)g(r,x-z)dz|^2dr\right)^{\frac{p}{2}}\\
&=:&\int_{Q_1\setminus Q_{12}}\int_{Q_1\setminus Q_{12}}(J_1+J_2+J_3)dtdxdsdy.
  \eess

Estimate of $J_1$. By using the H\"{o}lder continuous of $g$, i.e.,
   \bess
|g(r,x-z)-g(0,0)|&\leq& C_g\max\left\{r^{\frac{1}{2}},|x-z|\right\}^\beta\\
&\leq&C(g,\beta)(T^{\frac{\beta}{2}}+|x-x_1|^\beta+|x_1|^\beta+|z|^\beta)\\
&\leq&C(g,\beta)(T^{\frac{\beta}{2}}+c_1^\beta+|x_1|^\beta+|z|^\beta),
   \eess
and (\ref{3.8}), we have
   \bess
J_1&=&C(p)
\left(\int_0^s|\int_{\mathbb{R}^d}|K(t-r,z)-K(s-r,z)||g(r,x-z)|dz|^2
dr\right)^{\frac{p}{2}}\\
&\leq&C(p,\beta,T)
\left(\int_0^s|\int_{\mathbb{R}^d}|K(t-r,z)-K(s-r,z)|(1+|z|^\beta)dz|^2
dr\right)^{\frac{p}{2}}\\
&&+c_1^{\beta p} C(p,\beta)
\left(\int_0^s\int_{\mathbb{R}^d}|K(t-r,z)-K(s-r,z)|
dr\right)^{\frac{p}{2}}\\
&\leq&C(p,\beta,T)(1+c_1^{\beta p})(t-s)^{{\frac{\gamma_1p}{2}}}.
   \eess
Here and in the rest part of the proof, we write the constant depending
on $\|g\|_{ C^\beta(\mathbb{R}_+\times\mathbb{R}^d))}$ as $C(\beta)$ for simplicity.
The condition (\ref{3.9}) and
   \bess
|g(r,x-z)-g(r,y-z)|\leq C_g|x-y|^\beta
   \eess
imply the following derivation
   \bess
J_2&=&C(p)\int_Q\int_Q\left(\int_0^s|\int_{\mathbb{R}^d}|K(s-r,z)||g(r,x-z)-g(r,y-z)|dz|^2dr\right)^{\frac{p}{2}}\\
&\leq&C(p,g)\int_Q\int_Q\left(\int_0^s|\int_{\mathbb{R}^d}|K(r,z)||x-y|^\beta dz|^2dr\right)^{\frac{p}{2}}\\
&\leq& C(N_0,p,g,\beta)|x-y|^{\beta p}.
   \eess

Estimate of $I_3$. By using the property $g(0,0)=0$ and (\ref{3.10}), we get
   \bess
J_3&=&C(p)\left(\int_s^t|
\int_{\mathbb{R}^d}K(t-r,z)g(r,x-z)dz|^2dr\right)^{\frac{p}{2}}\\
&\leq& C\left(\int_s^t\Big|
\int_{\mathbb{R}^d}|K(r,z)|
(T+|x-x_1|^\beta+|x_1|^\beta+|z|^\beta) dz\Big|^2dr\right)^{\frac{p}{2}}\\
&\leq&C(p,T, \beta)
\left(\int_s^t\Big|\int_{\mathbb{R}^d}|K(t-r,z)|(1+|z|^\beta) dz
\Big|^2dr\right)^{\frac{p}{2}}\\
&&+C(p,T,\beta)|x-y|^{\beta p}
\left(\int_s^t\Big|\int_{\mathbb{R}^d}|K(t-r,z)| dz
\Big|^2dr\right)^{\frac{p}{2}}\\
&\leq&C(p,T, \beta) (t-s)^{\frac{\gamma_2p}{2}}(1+|x-y|^{\beta p}).
   \eess

Noting that $(t,x)\in Q_{1}$ and $(s,y)\in Q_{1}$, we have
   \bess
0\leq t-s\leq 2{c_1}^2\ \ \ {\rm and}\ \ \ |x-y|\leq|x-x_1|+|y-x_1|\leq2{c_1}.
   \eess
Using the above inequality and the properties of $A$-type domain, we deduce
  \bess
 \int_{Q_1\setminus Q_{12}}\int_{Q_1\setminus Q_{12}}J_1dtdxdsdy&\leq&C(p,T,\beta)(1+c_1^{\beta p}){c_1}^{\gamma_1p}|Q_1\setminus Q_{12}|^2;\\
\int_{Q_1\setminus Q_{12}}\int_{Q_1\setminus Q_{12}}J_2dtdxdsdy&\leq& C(C_0,p,g,\beta){c_1}^{\beta p}|Q_1\setminus Q_{12}|^2;\\
\int_{Q_1\setminus Q_{12}}\int_{Q_1\setminus Q_{12}}J_3dtdxdsdy&\leq&C(p,T,\beta)|Q_1\setminus Q_{12}|^2{c_1}^{\gamma_2p}(1+{c_1}^{\beta p}).
   \eess
Combining the estimates of $J_1,J_2$ and $J_3$, we get
 \bess
&&\mathbb{E}\int_{Q_1\setminus Q_{12}}\int_{Q_1\setminus Q_{12}}|u(t,x)-u(s,y)|^pdtdxdsdy\\
&\leq&C(\beta,C_0,T,p)|Q_1\setminus Q_{12}|^2(c_1^{\beta p}+1)(c_1^{\beta p}+c_1^{\gamma_1p}+c_1^{\gamma_2p}).
   \eess
Since $D$ is a $A$-type bounded domain, we have for $2c_1\leq diam D$,
   \bess
&&A|Q_{c_1}(t_1,x_1)|\leq |Q_1|\leq |Q_{c_1}(t_1,x_1)|\\
&&A|Q_{c_1}(t_1,x_1)\setminus Q_{c_1}(t_2,x_1)|\leq |Q_1\setminus Q_{12}|\leq |Q_{c_1}(t_1,x_1)\setminus Q_{c_1}(t_2,x_1)|.
   \eess
We remark that
   \bess
&&|Q_{c_1}(t_1,x_1)|\approx Cc_1^{d+2},\\
&&|Q_{c_1}(t_1,x_1)\setminus Q_{c_1}(t_2,x_1)|\leq Cc_1^{d}[c_1^2\wedge(t_1-t_2)],
   \eess
where $C$ is a positive constant which does not depend on
$c_1$. Noting that $Q_1\setminus Q_{12}\subset Q_1$ and taking $0<\delta<1$, we have
  \bess
&&\mathbb{E}
\int_{Q_1\setminus Q_{12}}\int_{Q_1\setminus Q_{12}}|\mathcal {K}g(t,x)-\mathcal {K}g(s,y)|^pdtdxdsdy\\
&\leq&C(\beta,C_0,D,d,T)|Q_1|^{2+\frac{\gamma p-2\delta}{d+2}}|t_1-t_2|^\delta,
   \eess
where $\gamma=\min\{\gamma_1,\gamma_2,\beta\}$.

Similarly, we can get
 \bess
&&\mathbb{E}\int_{Q_{12}\setminus Q_{1}}\int_{Q_{12}\setminus Q_{1}}|u(t,x)-u(s,y)|^pdtdxdsdy\\
&\leq&C(\beta,C_0,D,d,T)|Q_1|^{2+\frac{\gamma p-2\delta}{d+2}}|t_1-t_2|^\delta.
   \eess
Due to the fact that $I_{12}\leq0$, we have
   \bess
\mathbb{E}I_1\leq C(\beta,C_0,D,d,T)|t_1-t_2|^\delta,
  \eess
where $\theta=1+\frac{\gamma p-2\delta}{d+2}$.

Next, similar to the proof of Theorem \ref{t3.1}, one can estimate $I_2$ and $I_3$ as
followings
    \bess
&&\mathbb{E}I_2=\mathbb{E}[F(t_2,x_1,c_1)-F(t_2,x_2,c_1)]\leq C(\beta,C_0,D,d,T)|x_1-x_2|^\delta,\\
&&\mathbb{E}I_3=\mathbb{E}[F(t_2,x_2,c_1)-F(t_2,x_2,c_2)]\leq C(\beta,C_0,D,d,T)|c_1-c_2|^\delta,
   \eess
where $\theta=1+\frac{\gamma p-\delta}{d+2}$.

Therefore, we have
   \bess
 &&\mathbb{E}|F(t_1,x_1,c_1)-F(t_2,x_2,c_2)|^q\\
 &\leq& C(C_0,D,d,T)\|g\|_{ C^\beta(\mathbb{R}_+\times\mathbb{R}^d))}^{q}(|t_1-t_2|+|x_1-x_2|+|c_1-c_2|)^{\delta q},
    \eess
where $\theta=1+\frac{\beta p-2\delta}{d+2}$, $(t_i,x_i)\in D_T$ and $0<c_i\leq d(D)$, $i=1,2$.
The rest proof of this theorem is exactly similar to that of Theorem \ref{t3.1} and we omit it here.
The proof of Theorem \ref{t3.2} is complete. $\Box$

Next, we consider the following equation
   \bes
\frac{\partial}{\partial t}u(t,x)=\Delta^\alpha u(t,x)+g(t,x)\dot{W}(t,x),\ \ \ u|_{t=0}=0,
   \lbl{3.13}\ees
where $\Delta^\alpha=-(-\Delta)^\alpha$ and $W_t$ is a
standard Brownian motion on a filtered probability space
$(\Omega,\mathcal {F},\{\mathcal {F}_t\}_{t\ge0},\mathbb{P})$.
\begin{theo}\lbl{t3.3} Let $D$ be an $A$-type bounded domain in $\mathbb{R}^{d+1}$ such that $\bar D\subset \mathcal {O}_T$.
Suppose that $g\in L^\infty_{loc}(\mathbb{R}_+;L^p(\Omega\times\mathbb{R}^d))$
 is $\mathcal {F}_t$-adapted process.
Set $d=1$. Assume that $\frac{1}{2}<\alpha\leq1$, $p>\frac{2}{2\alpha-1}$. Let $\beta>0$ be sufficiently
small such that $p(2\alpha-2\beta-1)>2$. Then, there is a mild solution $u$ of (\ref{3.13}) and
$u\in \mathscr{L}^{p,\theta}((D_T;\delta);L^p(\Omega))\cap L^p(\Omega;C^{\beta}(D_T))$. Moreover, it holds that
   \bes
&&\|u\|_{\mathscr{L}^{p,\theta}((D_T;\delta);L^p(\Omega))}\leq C\|g\|_{ L^\infty([0,T];L^p(\Omega\times\mathbb{R}))},\lbl{3.14}\\
&&\|u\|_{C^{\beta}(D_T;L^p(\Omega))}\leq C\|g\|_{ L^\infty([0,T];L^p(\Omega\times\mathbb{R}))},\lbl{3.15}
  \ees
where $\theta=1+\frac{\beta p}{3}$. Moreover,  taking $0<\delta<\beta p/2$ and $q>3/\delta$,
we have for $0<r<q$
   \bes
\|u\|_{L^r(\Omega;C^{\beta^*}(D_T))}\leq C\|g\|_{ L^\infty([0,T];L^p(\Omega\times\mathbb{R}))},\lbl{3.16}
   \ees
where $\beta^*=\beta-2\delta/p$.
  \end{theo}

{\bf Proof.}  The existence of mild solution of (\ref{3.13}) is a classical result under
the above assumptions. Now we prove the inequality (\ref{3.14}).
Due to the definition of Companato space, it suffices to show that
   \bess
[u]_{\mathscr{L}^{p,\theta}((D_T;\delta);L^p(\Omega))}<\infty.
  \eess
Direct calculus shows that
    \bess
[u]^p_{\mathscr{L}^{p,\theta}((D_T;\delta);L^p(\Omega))}&\leq&
\sup_{D(X,c),X\in D_T,0<c\leq d(D)}\frac{1}{|D(X,c)|^{1+\theta}}\\
&&\times\mathbb{E}
\int_{D(X,c)}\int_{D(X,c)}|u(t,x)-u(s,y)|^pdtdxdsdy\\
&\leq&\sup_{D(X,c),X\in D_T,0<c\leq d(D)}\frac{1}{|D(X,c)|^{1+\theta}}\\
&&\times\mathbb{E}
\int_{D(X,c)}\int_{D(X,c)}\Big|\int_0^t\int_{\mathbb{R}}K(t-r,x-z)g(r,z)dzdW(r)\\
&&-\int_0^s\int_{\mathbb{R}}K(s-r,y-z)g(r,z)dW(dr,dz)\Big|^p\\
&:=&\sup_{D(X,c),X\in D_T,0<c\leq d(D)}\frac{1}{|D(X,c)|^{1+\theta}}
\int_{D(X,c)}\int_{D(X,c)}\mathbb{E}\Upsilon dtdxdsdy.
   \eess

Set $t\geq s$. We have the following estimates
    \bess
\mathbb{E}\Upsilon
&\leq&C\mathbb{E}\Big|\int_0^s\int_{\mathbb{R}}(K(t-r,x-z)-K(s-r,y-z))g(r,z)W(dr,dz)\Big|^p\\
&&
+C\mathbb{E}\Big|\int_s^t\int_{\mathbb{R}}K(t-r,x-z)g(r,z)W(dr,dz)\Big|^p\\
&\leq&C\mathbb{E}\Big|\int_0^s\int_{\mathbb{R}}(K(t-r,x-z)-K(s-r,y-z))^2g^2(r,z)dzdr\Big|^{\frac{p}{2}}\\
&&
+C\mathbb{E}\Big|\int_s^t\int_{\mathbb{R}}K^2(t-r,x-z)g^2(r,z)dzdr\Big|^{\frac{p}{2}}\\
&=:&C(H_1+H_2).
   \eess

{\bf Estimate of $H_1$}.

Take $\beta>0$ satisfying $(2\alpha-2\beta-1)p-2\geq0$.
By using the Proposition \ref{p2.3}, and H\"{o}lder inequality, we have
   \bess
H_1 &=&\mathbb{E}\Big|\int_0^s\int_{\mathbb{R}}(K(t-r,x-z)-K(s-r,y-z))^2g^2(r,z)dzdr\Big|^{\frac{p}{2}}\\
 &\leq&C\mathbb{E}\Big|\int_0^s\int_{\mathbb{R}}|K(t-r,x-z)-K(s-r,x-z)|^2\cdot|g^2(r,z)|dzdr\Big|^{\frac{p}{2}}\\
 &&+C\mathbb{E}\Big|\int_0^s\int_{\mathbb{R}}(K(s-r,x-z)-K(s-r,y-z))^2\cdot g^2(r,z) dzdr\Big|^{\frac{p}{2}}\\
 &=:&H_{11}+H_{12}.
    \eess
 For $H_{11}$, we have
   \bess
 H_{11}
 &\leq&C(t-s)^{\frac{\beta p}{2}}\mathbb{E}\Big|\int_0^s\left(\int_{\mathbb{R}}
 |\frac{\partial^{\frac{\beta }{2}} K}{\partial t^{\frac{\beta }{2}}}(\xi-r,x-z)|^qdz \right)^{\frac{2}{q}}
\|g(r)\|^2_{L^p(\mathbb{R})}dr\Big|^{\frac{p}{2}}\\
&\leq&C(t-s)^{\frac{\beta p}{2}}\|g\|_{L^p(\Omega;L^\infty([0,T];L^p(\mathbb{R})))}^p\left[\int_0^s
\left(\int_{\mathbb{R}}|\frac{\partial^{\frac{\beta}{2}} K}{\partial t^{\frac{\beta}{2}}}(\xi-r,x-z)|^qdz \right)^{\frac{2}{q}}dr\right]^{\frac{p}{2}},
   \eess
where $q=2p/(p-2)$, $\xi=\theta t+(1-\theta)s$, $0<\theta<1$ and we used the
following fact
   \bess
&&\int_0^s
\left(\int_{\mathbb{R}}|\frac{\partial^{\frac{\beta}{2}} K}{\partial t^{\frac{\beta}{2}}}(\xi-r,x-z)|^qdz \right)^{\frac{2}{q}}dr
\\
&\leq&C\int_0^s\left(\int_0^{(\xi-r)^{\frac{1}{2\alpha}}}(\xi-r)^{-\frac{q+2q\alpha\beta}{2\alpha}}|z|^{q\alpha\beta}d|z|\right.\\
&&\left.
+\int_{(\xi-r)^{\frac{1}{2\alpha}}}^\infty(\xi-r)^q|z|^{-(q+2q\alpha+2q\alpha\beta)}|z|^{q\alpha\beta}d|z| \right)^{\frac{2}{q}}dr\\
&\leq&C\left[(\theta(t-s))^{\frac{1-q+q\alpha(1-\beta)}{q\alpha}}+\xi^{\frac{1-q+q\alpha(1-\beta)}{q\alpha}}\right]\\
&\leq&C
  \eess
because using $q=2p/(p-2)$, we have
   \bess
1-q+q\alpha(1-\beta)>0\Leftrightarrow p(2\alpha-2\alpha\beta-1)>2\Leftarrow p(2\alpha-2\beta-1)>2.
  \eess

 For $H_{12}$, by using the fractional mean value formula again,
we have
  \bess
 H_{12}&\leq&
C|x-y|^{\beta p}\|g\|_{L^p(\Omega;L^\infty([0,T];L^p(\mathbb{R})))}^p\Big|\int_0^s \left(\int_{\mathbb{R}}[K^{(\beta)}(s-r,\xi-z)]^{q}dzdr\right)^{\frac{2}{q}}
dr\Big|^{\frac{p}{2}}\\
&\leq&C|x-y|^{\beta p}\|g\|_{L^p(\Omega;L^\infty([0,T];L^p(\mathbb{R})))}^p\left[\int_0^s(s-r)^{-\frac{d(q-1)+\beta q}{q\alpha}}dr\right]^{\frac{p}{2}}\\
&\leq&C|x-y|^{\beta p},
   \eess
where $q=2p/(p-2)$, $\xi=\theta x+(1-\theta)y$ and we used the following inequality
  \bess
\int_0^s(s-r)^{-\frac{d(q-1)+\beta q}{q\alpha}}dr=\frac{q\alpha}{q(\alpha-\beta)-(q-1)}s^{\frac{q(\alpha-\beta)-(q-1)}{q\alpha}}
\leq C
   \eess
provided that $(2\alpha-2\beta-1)p-2\geq0$.

{\bf Estimate of $H_2$}.

Similar to the estimate of $H_1$, we have
   \bess
H_2&=&\mathbb{E}\Big|\int_s^t\int_{\mathbb{R}^d}K^2(t-r,x-z)g^2(r,z)dzdr\Big|^{\frac{p}{2}}\\
&\leq&\|g\|_{L^p(\Omega;L^\infty([0,T];L^p(\mathbb{R}^d)))}^p
\left[\int_s^t\left(\int_{\mathbb{R}^n}|K(t-r,x-z)|^qdz \right)^{\frac{2}{q}}dr\right]^{\frac{p}{2}}\\
&\leq&C\|g\|_{L^p(\Omega;L^\infty([0,T];L^p(\mathbb{R}^d)))}^p(t-s)^{\frac{q\alpha-(q-1)d}{q\alpha}\times \frac{p}{2}}
  \eess
provided that $p(2\alpha-1)> 2$. Indeed, by using $q=\frac{2p}{p-2}$, we have
   \bess
q\alpha-(q-1)>0\ \ \Longleftrightarrow \ \ p(2\alpha-1)> 2.
  \eess

Combining the assumption of $p$, we have
    \bess
H_2\leq C(t-s)^{\frac{p(2\alpha-1)-2}{2\alpha}}.
   \eess
The rest proof is similar to that of \ref{t3.1} and we omit it here.  $\Box$

\section{H\"{o}lder estimate on a bounded domain}\setcounter{equation}{0}
In this section, we consider the SPDEs of the following form
   \bes\left\{
\begin{array}{llll}
du=Audt+g(t,x)dW_t, \ \ \ (t,x)\in (0,\infty)\times D,\\
u|_{\partial D}=0,\\
u_{t=0}=0,
  \end{array}\right.\lbl{5.1}\ees
where $D$ is a smooth bounded domain in $\mathbb{R}^d$, $W_t$ is
standard one-dimensional Brownian motion, and $g$ is progressively measurable
$L^\infty$- or $L^p$-function.

Throughout this section, we assume that $A$ is a uniformly elliptic second-order
differential operator of the form
   \bess
A=a_{ij}\frac{\partial}{\partial x_i}\frac{\partial}{\partial x_j}+b_i(x)\frac{\partial}{\partial x_i}+c(x)
   \eess
with smooth coefficients. Furthermore, we assume that at least one of the
following two assumptions holds:
   \bess
&&B^\infty:\ \ \ \|g\|_{L^\infty([0,T];L^p(\Omega;L^\infty(D)))}<\infty,\\
&&B^p:\ \ \ \ \|g\|_{L^\infty([0,T];L^p(\Omega\times D)))}<\infty.
   \eess
In order to obtain the H\"{o}lder estimate, we need the following Lemma.
Consider the following initial-boundary problem:
   \bes
\frac{\partial v}{\partial t}-Av=0,\ \ v|_{t=0}=F(x),\ \ v|_{\partial D}=0,
   \lbl{5.2}\ees
and denote by $S_t$ the corresponding semigroup:
    \bess
v(t,\cdot)=(S_tF)(\cdot),\ \ \ F=F(\cdot).
   \eess

  \begin{lem}\lbl{l5.1} \cite[Lemma 1]{KNP} Let $|F(x)|<M$. Then, for any $\theta<1$,
the following estimates hold with $c>0$:
    \bess
&&\|v(t,\cdot)\|_{C^\theta(D)}\leq c(\theta)Mt^{-\theta/2}\exp(-ct),\\
&&|v(t+\delta,x)-v(t,x)|\leq c(\theta)Mt^{-\theta}\delta^\theta\exp(-ct).
   \eess
Moreover, if $\|F\|_{L^p(D)}\leq M$ and $p>1$, then
   \bess
&&\|v(t,\cdot)\|_{C^\theta(D)}\leq c(\theta)Mt^{-\theta/2-d/(2p)}\exp(-ct),\\
&&|v(t+\delta,x)-v(t,x)|\leq c(\theta)Mt^{-\theta-d/(2p)}\delta^\theta\exp(-ct).
   \eess
  \end{lem}

\begin{theo}\lbl{t5.1} Let $D_T$ be an $A$-type bounded domain in $\mathbb{R}^{d+1}$.

(i) Suppose that $B^p$ holds for $p>d$
 and that $0<\beta<1$ satisfies $(1-\beta)p-d\geq0$. Then, there is a mild solution $u$ of (\ref{5.1}) and
$u\in \mathscr{L}^{p,\theta}((D_T;\delta);L^p(\Omega))\cap L^p(\Omega;C^{\beta}(D_T))$. Moreover, it holds that
   \bess
&&\|u\|_{\mathscr{L}^{p,\theta}((D_T;\delta);L^p(\Omega))}\leq C\|g\|_{ L^\infty([0,T];L^p(\Omega\times D))},\\
&&\|u\|_{C^{\beta}(D_T;L^p(\Omega))}\leq C\|g\|_{ L^\infty([0,T];L^p(\Omega\times D))},
  \eess
where $\theta=1+\frac{\beta p}{d+2}$. Moreover,  taking $0<\delta<\beta p/2$ and $q>(d+2)/\delta$,
we have for $0<r<q$
   \bess
\|u\|_{L^r(\Omega;C^{\beta^*}(D_T))}\leq C\|g\|_{ L^\infty([0,T];L^p(\Omega\times\mathbb{R}^d))},
   \eess
where $\beta^*=\beta-2\delta/p$.

(ii) Suppose that $B^\infty$ holds for $p>1$. Then, there is a mild solution $u$ of (\ref{5.1}) and
$u\in \mathscr{L}^{p,\theta}((D_T;\delta);L^p(\Omega))\cap L^p(\Omega;C^{\beta}(D_T))$. Moreover, it holds that
   \bess
&&\|u\|_{\mathscr{L}^{p,\theta}((D_T;\delta);L^p(\Omega))}\leq C\|g\|_{ L^\infty([0,T];L^p(\Omega\times D))},\\
&&\|u\|_{C^{\beta}(D_T;L^p(\Omega))}\leq C\|g\|_{ L^\infty([0,T];L^p(\Omega\times D))},
  \eess
where $\theta=1+\frac{ p}{d+2}$. Moreover,  taking $0<\delta< p/2$ and $q>(d+2)/\delta$,
we have for $0<r<q$
   \bess
\|u\|_{L^r(\Omega;C^{\beta^*}(D_T))}\leq C\|g\|_{ L^\infty([0,T];L^p(\Omega\times\mathbb{R}^d))},
   \eess
where $\beta^*=1-2\delta/p$.
  \end{theo}

{\bf Proof.} The proof of this Theorem is exactly similar to that of Theorem \ref{t3.1} by using
Lemma \ref{l5.1}. We omit it to the readers. The proof is complete. $\Box$

\begin{remark}\lbl{r5.1} Theorem \ref{t5.1} does not hold for the nonlocal operator because
we did not have the similar properties of kernel function on bounded domain.

Comparing Theorem \ref{t5.1} with \cite[Theorems 1 and 2]{KNP}, we find the index of \cite{KNP}
is $\beta<\frac{1}{2}-\frac{d}{2p}$ for the case $B^p$ and the index in this paper is larger than that of
\cite{KNP}.
\end{remark}

\section{Applications and further discussions}\setcounter{equation}{0}

We first give an example for Theorem \ref{t3.2}.
Consider the equation (\ref{3.1}).
In our paper \cite{LGWW1}, by using
Proposition \ref{p2.3}, we got
the following result.
\begin{lem}\lbl{l4.1}Let $0\leq\epsilon< \alpha $. The following estimates hold.
   \bess
&&\int_0^s\left(\int_{\mathbb{R}^d}|\nabla^\epsilon p(t-r,z)-\nabla^\epsilon p(s-r,z)|(1+|z|^\beta)dz\right)^2dr\leq N(T,\beta)(t-s)^{\gamma},\\
&&\int_0^s\left(\int_{\mathbb{R}^d}|\nabla^\epsilon p(s-r,z)|dz\right)^2dr\leq N_0,\\
&& \int_s^t\left(\int_{\mathbb{R}^d}|\nabla^\epsilon p(t-r,z)|(1+|z|^\beta) dz\right)^2dr\leq N(T,\beta)(t-s)^{\gamma},
   \eess
where $\gamma=\frac{\alpha-\epsilon}{\alpha}$.
\end{lem}
Then applying Theorem \ref{t3.2}, we have the following result.

 \begin{theo}\lbl{t4.1} Let $0\leq\epsilon< \alpha $ and
$D_T$ be an $A$-type bounded domain in $\mathbb{R}^{d+1}$ such that
$\bar D_T\subset \mathcal {O}_T$.
Suppose that $g\in C^\beta(\mathbb{R}_+\times\mathbb{R}^d)$, $0<\beta<1$, is a non-random function
  and $g(0,0)=0$.
Then we have, for
$p\geq1$ and $\beta<\gamma$,
    \bess
&&\|\nabla^\epsilon u\|_{\mathscr{L}^{p,\theta}((D_T;\delta);L^p(\Omega))}\leq
C\|g\|_{ C^\beta(\mathbb{R}_+\times\mathbb{R}^d))},\nonumber\\
&&\|\nabla^\epsilon u\|_{C^{\beta}(D_T;L^p(\Omega))}\leq C\|g\|_{ C^\beta(\mathbb{R}_+\times\mathbb{R}^d))},
  \eess
where $\theta=1+\frac{\gamma p}{d+2}$ and $\gamma=\frac{\alpha-\epsilon}{\alpha}$.
Moreover,  taking $0<\delta<\gamma p/2$ and $q>(d+2)/\delta$,
we have for $0<r<q$
   \bess
\|\nabla^\epsilon u\|_{L^r(\Omega;C^{\beta^*}(D_T))}\leq C\|g\|_{ C^\beta(\mathbb{R}_+\times\mathbb{R}^d))},
   \eess
where $\beta^*=\gamma-2\delta/p$.
\end{theo}
In fact, one can use the factorization method to obtain the H\"{o}lder estimates of solutions
to the following equation
  \bess
du_t=[\Delta^\alpha u+f(t,x,u)]dt+g(t,x)dW_t,\ \ \ u|_{t=0}=u_0(x),
   \eess
where $\Delta^\alpha=-(-\Delta)^\alpha$, $\alpha\in(0,1]$ and $W_t$ is a
standard Brownian motion on a filtered probability space
$(\Omega,\mathcal {F},\mathcal {F}_t,\mathbb{P})$. About the factorization method,
see \cite{DMH}.

In addition, one can use the Kunita's first inequality to deal with a general case.
Let $(\Omega,\mathcal {F},\mathbb{F},\mathbb{P})$ be a complete probability space such that
$\{\mathcal {F}_t\}_{t\in[0,T]}$ is a filtration on $\Omega$ containing
all $P$-null subsets of $\Omega$ and $\mathbb{F}$ be the predictable $\sigma$-algebra associated with the filtration $\{\mathcal {F}_t\}_{t\in[0,T]}$. We are given a $\sigma$-finite measure space $(Z,\mathcal {Z},\nu)$ and a Poisson
random measure $\mu$ on $[0,T]\times Z$, defined on the stochastic basis. The compensator of
$\mu$ is Leb$\otimes\nu$, and the compensated martingale measure $\tilde N:=\mu-Leb\otimes\nu$.
The method used here is also suitable to the case that
    \bes
\mathcal {G}g(t,x)&=&\int_0^t\int_ZK(t,s,\cdot)\ast g(s,\cdot,z)(x)\tilde N(dz,ds)\nm\\
&=&\int_0^t\int_Z\int_{\mathbb{R}^d}K(t-s,x-y)g(s,y,z)dy\tilde N(dz,ds)
    \lbl{4.1}\ees
for $\mathbb{F}$-predictable processes $g:[0,T]\times\mathbb{R}^d\times Z\times\Omega\to\mathbb{R}$.

In the end of this section, we give a new criteria based on the following Proposition.
\begin{prop} \cite[Theorem 2.1]{RY} Let $\{X_t, t\in [0,1]\}$ be a Banach-valued stochastic field for which there exist three strictly positive constants $\gamma,c, \varepsilon$ such that
\bess
E[\sup_{0\leq t\leq 1}|X_t(x)-X_t(y)|^\gamma]\leq c|x-y|^{d+\varepsilon},
\eess
then there is a modification $\tilde{X}$ of $X$ such that
\bess
E\Big[\Big(\sup_{s\neq t}\frac{|\tilde{X}_t-\tilde{X}_s|}{|t-s|^\alpha}\Big)^\gamma\Big]<\infty
\eess
for every $\alpha\in [0,\varepsilon/\gamma)$. In particular, the paths of $\tilde{X}$ are H\"{o}lder continuous in $x$ of order $\alpha$.
\end{prop}

For applications, we need prove the Kolmogorov criterion with the following form.
\begin{theo} \label{t4.2} Let $\{X_t(x), x\in [0,1]^d, t\in [0,1]\}$ be a Banach-valued stochastic field for which there exist three strictly positive constants $\gamma,c, \varepsilon$ such that
\bess
E[\sup_{0\leq t\leq 1}|X_t(x)-X_t(y)|^\gamma]\leq c|x-y|^{d+\varepsilon},
\eess
then there is a modification $\tilde{X}$ of $X$ such that
\bess
E\Big[\sup_{0\leq t\leq 1}\Big(\sup_{x\neq y}\frac{|\tilde{X}_t(x)-\tilde{X}_t(y)|}{|x-y|^\alpha}\Big)^\gamma\Big]<\infty
\eess
for every $\alpha\in [0,\varepsilon/\gamma)$. In particular, the paths of $\tilde{X}$ are H\"{o}lder continuous in $x$ of order $\alpha$.
\end{theo}
\vskip2mm\noindent
\textbf{Proof.} Let $D_m$ be the set of points in $[0,1]^d$ whose components are equal to $2^{-m}i$ for some integral $i\in [0,2^m]$. The set $D=\cup_{m}D_m$ is the set of dyadic numbers. Let further $\Delta_m$ be the set of pairs $(x,y)$ in $D_m$ such that $|x-y|=2^{-m}$. There are $2^{(m+1)d}$ such pairs in $\Delta_m$.

Let us finally set $K_i(t)=\sup_{(x,y)\in \Delta_i}|X_t(x)-X_t(y)|$. The hypothesis entails that for a constant $J$,
\bess
E[\sup_{0\leq t\leq 1}K_i(t)^\gamma ]\leq \sum_{(x,y)\in \Delta_i}E[\sup_{0\leq t\leq 1}|X_t(x)-X_t(y)|^\gamma]\leq c2^{(i+1)d}2^{-i(d+\varepsilon)}=J2^{-i\varepsilon}.
\eess

For a point $x$ (resp. $y$) in $D$, there is an increasing sequences $\{x_m\}$ (resp. $\{y_m\}$) of points in $D$ such that $x_m$ (resp. $y_m$) is in $D_m$ for each $m$, $x_m\leq x$ ($y_m\leq y$) and $x_m=x$ ($y_m=y$) from some $m$ on. If $|x-y|\leq 2^{-m}$, then either $x_m=y_m$ or $(x_m,y_m)\in \Delta_m$ and in any case
\bess
X_t(x)-X_t(y)=\sum_{i=m}^\infty(X_t(x_{i+1})-X_t(x_{i}))+X_t(x_m)-X_t(y_m)-
\sum_{i=m}^\infty(X_t(y_{i+1})-X_t(y_{i})),
\eess
where the series are actually finite sums. It follows that
\bess
|X_t(x)-X_t(y)|\leq K_m+2\sum_{i=m+1}^\infty K_i(t)\leq 2\sum_{i=m}^\infty K_i(t).
\eess
As a result, setting $M_\alpha(t)=\sup\{ |X_t(x)-X_t(y)|/|x-y|^\alpha, \ x,y\in D, \ x\neq y\}$, we have
\bess
M_\alpha(t) &\leq&\sup_{m\in N}\Big\{ 2^{m\alpha}\sup_{|x-y|\leq 2^{-m}} |X_t(x)-X_t(y)|, \ x,y\in D, \ x\neq y\Big\}
\cr\cr&\leq & \sup_{m\in N}\Big\{ 2^{m\alpha+1}\sum_{i=m}^\infty K_i(t) \Big\}\cr\cr&\leq & 2 \sum_{i=0}^\infty2^{i\alpha} K_i(t).
\eess

For $\gamma\geq 1$ and $\alpha <\varepsilon/\gamma$, we get with $J^\prime=2J$,
\bess
[ E \sup_{0\leq t\leq 1}M_\alpha(t)^\gamma]^{1/\gamma} \leq
2 \sum_{i=0}^\infty2^{i\alpha} [E \sup_{0\leq t\leq 1}K_i(t)^\gamma]^{1/\gamma}
\leq J^\prime \sum_{i=0}^\infty2^{i(\alpha-\varepsilon/\gamma)} <\infty.
\eess
For $\gamma< 1$, the same reasoning applies to $[E \sup_{0\leq t\leq 1}M_\alpha(t)^\gamma]$ instead of $[E \sup_{0\leq t\leq 1}M_\alpha(t)^\gamma]^{1/\gamma}$.

It follows in particular that for almost every $\omega$, $X_t(\cdot)$ is uniformly continuous on $D$ and it is uniformly in $t$, so it make sense to set
\bess
\tilde{X}_t(x,\omega)=\lim_{y\in D, y\rightarrow x}X_t(y,\omega).
\eess
By Fatou's lemma and the hypothesis, $\tilde{X}_t(x)=X_t(x)$ a.s. and $\tilde{X}$ is clearly the desired modification. $\Box$

It is easy to see that one can use Theorem \ref{t4.2} to consider the equation (\ref{3.1}) and
(\ref{4.1})

\medskip

\noindent {\bf Acknowledgment} The first author was supported in part
by NSFC of China grants 11771123.

 \end{document}